\documentclass[12pt]{article}

\usepackage{fullpage}
\usepackage{amsmath} 
\usepackage{amssymb}
\usepackage{amstext}
\usepackage{amsthm}
\usepackage{amsfonts} 
\usepackage{mathtools}
\usepackage{enumerate}  
\usepackage{dsfont}
\usepackage[utf8]{inputenc}

\usepackage[pdftex]{graphicx} 
\usepackage{hyperref}

\theoremstyle{plain}

\newtheorem{theorem}{Theorem}
\newtheorem*{theorem*}{Theorem}
\newtheorem{lemma}[theorem]{Lemma}
\newtheorem{corollary}[theorem]{Corollary}
\newtheorem{claim}[theorem]{Claim}

\theoremstyle{definition}
\newtheorem{definition}{Definition}

\newcommand{\Q}{\mathbb Q}

\DeclareMathOperator{\OA}{OA}

\newcommand\ints{{\mathbb Z}}
\newcommand\re{{\mathbb R}}
\newcommand\rats{{\mathbb Q}}
\newcommand\Zw{{\mathbf w}}

\newcommand{\pmat}[1]{\begin{pmatrix}#1\end{pmatrix}}

\newcommand\comp[1]{{\mkern2mu\overline{\mkern-2mu#1}}}
\newcommand\diff{\mathbin{\mkern-1.5mu\setminus \mkern-1.5mu}}

\begin{document}

\title{State transfer in strongly regular graphs with an edge perturbation}

\author{Chris Godsil\footnote{University of Waterloo, Waterloo, Canada. email: \url{cgodsil@uwaterloo.ca}. C. Godsil gratefully acknowledges the support of the Natural Sciences and Engineering Council of Canada (NSERC), Grant No. RGPIN-9439.}, Krystal Guo\footnote{Universit\'{e} libre de Bruxelles, Brussels, Belgium. Most of this research was done when K. Guo was a post-doctoral fellow at University of Waterloo. email: \url{guo.krystal@gmail.com}}, Mark Kempton\footnote{Center of Mathematical Sciences and Applications, Harvard University, Cambridge MA, USA. Supported by Air Force Office of Scientific Research Grant FA9550-13-1-0097: Geometry and Topology of Complex Networks. email: \url{mkempton@cmsa.fas.harvard.edu}}, Gabor Lippner\footnote{Northeastern University, Boston, MA, USA. email: \url{g.lippner@neu.edu}}}

\date{  \today  }
\maketitle
\begin{abstract} Quantum walks, an important tool in quantum computing, have been very successfully investigated using techniques in algebraic graph theory. 
We are motivated by the study of state transfer in continuous-time quantum walks, which is understood to be a rare and interesting phenomenon. We consider a perturbation on an edge $uv$ of a graph where we add a weight $\beta$ to the edge and a loop of weight $\gamma$ to each of $u$ and $v$. We characterize when for this perturbation results in strongly cospectral vertices $u$ and $v$. Applying this to strongly regular graphs, we give infinite families of strongly regular graphs where some perturbation results in perfect state transfer. Further, we show that, for every strongly regular graph, there is some perturbation which results in pretty good state transfer. We also show for any strongly regular graph $X$ and edge $e \in E(X)$, that $\phi(X\diff  e)$ does not depend on the choice of $e$. \end{abstract}

\section{Introduction and preliminary definitions}

A quantum walk is a quantum analogue of a random walk on a graph. Childs shows in \cite{ChildsUniversalQComputation} that the continuous-time quantum walk is an universal computational primitive. The transition matrix of a continuous-time quantum walk on a graph $X$ is a matrix-valued function in time, $U(t)$, given as follows:
\[ U(t) = \exp(itA)
\]
where $A$ is the adjacency matrix of $X$. 

Quantum walks on graphs have proven to be highly valuable tools in quantum computation and information theory \cite{FarhiGutmann}.  They can be used to describe the fidelity of information transfer in networks of interacting qubits. A natural question, originally proposed by Bose \cite{BoseQuantumComPaths}, is whether such information transfer between two nodes in a graph can be perfect.

\begin{definition}
Let $X$ be a graph with vertices $u$ and $v$.
\begin{enumerate}[(i)]
\item We say that $X$ admits \emph{perfect state transfer} from $u$ to $v$ if there is some time $t>0$ such that
\[
|U(t)_{u,v}| = 1.
\]
\item We say  $X$ admits \emph{pretty good state transfer} from $u$ to $v$ if, for any $\epsilon>0$, there is a time $t>0$ such that
\[
|U(t)_{u,v}| > 1-\epsilon.
\]
\end{enumerate}
\end{definition} 

Throughout this paper, we will sometimes use the abbreviation pst to refer to perfect state transfer, and pgst to refer to pretty good state transfer.

Perfect state transfer has been studied in paths \cite{ChristandlPSTQuantumSpinNet2,GodsilStateTransfer12}, circulants \cite{BasicCirculant}, cubelike graphs \cite{GodsilCheungPSTCubelike}, and various other infinite families of graphs. See \cite{KendonTamon11} for a survey. It is classified for all known distance regular graphs \cite{CoutinhoGodsilGuoVanhove2}, where perfect state transfer in strongly regular graphs is also classified. Perfect state transfer is understood to be a rare phenomenon (see \cite{GodsilPerfectStateTransfer12}) and new examples are of interest. 

Pretty good state transfer has also been studied in several classes of graphs. Results in this area tend to use number-theoretic techniques and results about almost periodic functions. In \cite{GodsilKirklandSeveriniSmithPGST}, the authors show that $P_n$ admits pretty good state transfer between end vertices if and only if  $n = 2^m - 1$ where $m$ is an integer, or $n = p - 1$ or $2p - 1$, where $p$ is a prime. Pretty good state transfer has also been studied in \cite{BanchiCoutinhoGodsil, VinetZhedanovAlmost, CouGuovBo2017}. 

For a graph $X$, we will denote the characteristic polynomial of its adjacency matrix $A(X)$ by $\phi(X,t)$.  We say that vertices $u$ and $v$ of $X$ are \emph{cospectral} if $\phi(X\diff  u,t) = \phi(X\diff  v,t)$.   If $E_r$ are the idempotents in the spectral decomposition for $A(X)$, then it is known that $u$ and $v$ are copectral if and only if $(E_r)_{u,u} = (E_r)_{v,v}$ for all $r$. Vertices $u$ and $v$ are said to be \emph{strongly cospectral} if $E_re_u = \pm E_re_v$ for all $r$.  It is well-known (see, for instance \cite{GodsilPerfectStateTransfer12,BanchiCoutinhoGodsil}) that strong cospectrality is a necessary condition for both perfect and pretty good state transfer between two vertices.  

A \textsl{strongly regular graph} is a graph $X$ which is neither empty nor complete, such that every vertex has degree $k$, every pair of adjacent vertices has $a$ common neighbours and every pair of non-adjacent vertices has $c$ common neighbours. If $n$ is the number of vertices of the graph, then the tuple $(n,k,a,c)$ is said to be the \textsl{parameter set} of the strongly regular graph. We will frequently use the abbreviation SRG to refer to a strongly regular graph.  One of the unique properties of strongly regular graphs is that the adjacency matrix of a strongly regular graph has only 3 distinct eigenvalues. Throughout, we will be denoting these eigenvalues as by $k$, $\theta$, and $\tau$.  Here $k$ is a simple eigenvalue, equal to the degree of each vertex. We refer to \cite{GodsilRoyle} for further background on strongly regular graphs. 

Strongly regular graphs are walk-regular, which means that every pair of vertices is cospectral.  This makes them good candidates for the study of perfect state transfer.  However, because of high multiplicities of eigenvalues, these pairs of vertices are typically not strongly cospectral. In \cite{CoutinhoGodsilGuoVanhove2}, the authors show that a strongly regular graph $X$ admits perfect state transfer if and only if $X$ is isomorphic to the complement of the disjoint union of an even number of copies of $K_2$.  In this paper, we will study perturbations of strongly regular graphs.  We examine when these perturbations yield strongly cospectral pairs of vertices, and further, when these perturbations yield prefect or pretty good state transfer.  

Given a pair of vertices $u$ and $v$, we will add a weighted edge between $u$ and $v$ with weight $\beta$, and will put weighted loops on the vertices, with weight $\gamma$.  For a graph $X$, we will denote by $X^{\beta,\gamma}$ the graph perturbed in this way.  

The weighted loops, corresponding to perturbing the diagonal of the adjacency matrix, have also been called and \emph{energy potential} on the vertices.  Perfect and pretty good state transfer in the presence of a potential have been studied in \cite{KemptonLippnerYauPath,KemptonLippnerYauInvol}.

In this paper, we give an expression for the characteristic polynomial of $X^{\beta,\gamma}$ in Theorem \ref{thm:beta}  and determine when $u$ and $v$ are strongly cospectral in $X^{\beta,\gamma}$ in Theorem \ref{thm:perturb-str-cospec}, for any graph $X$. For strongly regular graphs, we are able to give sufficient conditions on $X$, $\beta$ and $\gamma$ for perfect state transfer to occur in $X^{\beta,\gamma}$ from $u$ to $v$. In particular, we prove the following results:

\begin{theorem}\label{mainthm:pst}
Let $X$ be a strongly regular graph with eigenvalues $k,\theta,\tau$.  Then if $k\equiv\theta\equiv\tau \pmod{4}$ is odd (resp. even), then for any pair of non-adjacent (resp. adjacent) vertices $u,v$ of $X$ there is a choice of $\beta,\gamma$ such that there is perfect state transfer between $u$ and $v$ in $X^{\beta,\gamma}$.  
\end{theorem}

Moreover, we show that there are infinite families of strongly regular graphs $X$ satisfying $k\equiv\theta\equiv\tau \pmod{4}$, giving infinitely many examples where perfect state transfer occurs in $X^{\beta,\gamma}$.

\begin{theorem}\label{mainthm:pgst}
Let $X$ be a strongly regular graph.  Then for any pair of vertices $u,v$ of $X$, there is a choice of $\beta,\gamma$ such that there is pretty good state transfer between $u$ and $v$ in $X^{\beta,\gamma}$.
\end{theorem}

Further, the choice of $\beta,\gamma$ is dense in the reals.

The remainder of the paper is organized as follows.  In Section \ref{sec:walkgen} we give some technical preliminaries regarding walk generating functions.  Section \ref{sec:perturb} investigates how the perturbation on a pair of vertices affects the characteristic polynomial, and characterizes when this perturbation yields strongly cospectral vertices.  In Section \ref{sec:srg}, we apply these tools specifically to strongly regular graphs, giving explicit expressions for the characteristic polynomial of $X^{\beta,\gamma}$ when $X$ is strongly regular.  In Section \ref{sec:pst}, we prove our main result, Theorem \ref{mainthm:pst}, on perfect state transfer.  Section \ref{sec:example} contains examples of our main theorem, and shows there are infinitely many strongly regular graphs satisfying the condition of Theorem \ref{mainthm:pst}.  Finally, in Section \ref{sec:pgst} we prove Theorem \ref{mainthm:pgst}.

\section{Walk generating functions}\label{sec:walkgen}

Let $X$ be a graph and $u,v$ be vertices of $X$. The \textsl{walk generating function} from $u$ to $v$ in $X$, denoted $W_{u,v}(X,t)$, is the generating function where the coefficient of $t^k$ is the number of walks from $u$ to $v$ in $X$. In other words,
\[
[t^k] W_{u,v}(X,t) = (A^k)_{u,v}
\]
and so
\[
W_{u,v}(X,t) = (I + tA + t^2A^2 + \cdots)_{u,v} = (I-tA)^{-1}_{u,v}.
\]
Following standard techniques (see \cite{GodsilAlgebraicCombinatorics}), we obtain immediately that
\begin{equation}\label{eq:waa}
W_{u,u}(X,t) = \frac{t^{-1} \phi(X\diff  u, t^{-1})}{ \phi(X, t^{-1})}.
\end{equation}
To find an expression for the off-diagonal entries of $(I-tA)^{-1}$, we let $\phi_{u,v}(X,t)$ be the polynomial such that 
\[
(I-tA)^{-1}_{u,v} = W_{u,v}(X,t) =  \frac{t^{-1} \phi_{u,v}(X, t^{-1})}{ \phi(X, t^{-1})}.
\]
The following lemma is found in \cite{GodsilAlgebraicCombinatorics}.

\begin{lemma}\cite{GodsilAlgebraicCombinatorics}\label{lem:wab} $ \phi_{u,v}(X, t)^2 = \phi(X\diff  u,t)\phi(X\diff  v,t) - \phi(X,t)\phi(X\diff  \{u,v\},t).$  \end{lemma}

Let $\theta_1, \ldots, \theta_m$ be the distinct eigenvalues of $A(X)$ and supposed $A(X)$ has the following spectral decomposition:
\[
A(X) = \sum_{r=1}^m \theta E_{\theta}
\]
where $E_{\theta}$ is the eigenmatrix representing the idempotent projection onto the $\theta_r$ eigenspace. We can obtain the following expression for $\phi(X\diff  u,t)$:
\begin{equation}\label{eq:phiaEr}
\frac{\phi(X\diff  u,t)}{\phi(X,t)} = \sum_{r=1}^m \frac{(E_r)_{u,u}}{t-\theta_r}.
\end{equation}
Similarly,
\begin{equation}\label{eq:phiabEr}
\frac{\phi_{uv}(X,t)}{\phi(X,t)} = \sum_{r=1}^m \frac{(E_r)_{u,v}}{t-\theta_r}.
\end{equation}

\section{Eigenvalues of the perturbation with potential}\label{sec:perturb}

Let $X$ be a graph with vertices  $u$ and $v$. Let $\beta$ and $\gamma$ fixed real numbers and let $H = \beta(e_u e_v^T + e_ve_u^T) + \gamma(e_u e_u^T + e_v e_v^T)$. In other words, $H$ is the matrix whose principal submatrix indexed by $\{u,v\}$ is 
\[
C := \pmat{\gamma & \beta \\ \beta & \gamma}
\]
and whose other entries are all equal to zero. Let $X^{\beta,\gamma}_{u,v}$ be the weighted graph whose adjacency matrix is $A + H$. We will refer to $X^{\beta,\gamma}_{u,v}$ as the \textsl{perturbation} of $X$ at vertices $u,v$ with weights $\beta,\gamma$. When the choice of vertices $u$ and $v$ is clear, we will omit the subscript and write simply $X^{\beta,\gamma}$. 

\begin{theorem} \label{thm:beta} If $u$ and $v$ are distinct vertices of $X$, then
\[\phi(X^{\beta,\gamma},t) 
= \phi(X,t) -2 \beta\phi_{uv} - \gamma (\phi_{v}+\phi_{u})
 +  (\gamma^2 -\beta^2) \phi(X\diff \{u,v\},t) .
\]
Further, if $u$ and $v$ are cospectral, then 
\[
\frac{\phi(X^{\beta,\gamma},t)}{ \phi(X,t)} = \left( 1 -\frac{(\beta-\gamma)(\phi_{uv}-\phi_{u})}{\phi(X,t)}\right) \left(1 - \frac{(\beta+\gamma)(\phi_{uv} + \phi_{u}) }{\phi(X,t)} \right). 
\]
\end{theorem} 

\proof
Let $E = \pmat{e_u & e_v}$. We see that $H = ECE^T$. We obtain that the following:
\[
\begin{split}
\phi(X^{\beta,\gamma},t) &= \det(tI - A - H)\\
&=  \det((tI - A)(I -(tI - A)^{-1} H)) \\
&= \phi(X,t) \det((I -(tI - A)^{-1} H)) \\
&= \phi(X,t) \det((I -(tI - A)^{-1} ECE^T)). 
\end{split}
\]
Note, for matrices $M$ and $N$ of appropriate orders, we have that 
\[\det(I-tNM) = \det(I-tMN).\]
We apply apply this property of determinants to obtain that 
\[ \det{(I -(tI - A)^{-1}  ECE^T)} = \det{(I - CE^T(tI - A)^{-1} E)}.
\]
For convenience, we will abbreviate $\phi(X\diff  u,t)$, $\phi(X\diff  v,t)$ and $\phi_{u,v}(X, t)$, as  $\phi_u$, $\phi_v$ and $\phi_{uv}$, respectively. We have 
\[
\begin{split}
CE^T(tI - A)^{-1} E &= C \pmat{e_u^T \\ e_v^T} (tI - A)^{-1} \pmat{e_u & e_v} \\
&=  \frac{1}{\phi(X,t)}  \pmat{\gamma & \beta \\ \beta & \gamma} \pmat{\phi_u & \phi_{uv} \\ \phi_{uv} & \phi_v}\\
&= \frac{1}{\phi(X,t)} \pmat{\beta\phi_{uv} + \gamma \phi_{u} & \beta\phi_v + \gamma  \phi_{uv} \\ \beta\phi_u +  \gamma\phi_{uv} & \beta \phi_{uv} +  \gamma\phi_{v}} ,
\end{split}
\]  
where we applied (\eqref{eq:waa}) to obtain the second line.
We obtain that
\[\det{(I - CE^T(tI - A)^{-1} E)} = \det \pmat{1 - \frac{1}{\phi(X,t)} (\beta\phi_{uv} + \gamma \phi_{u}) & \frac{- 1}{\phi(X,t)}(\beta\phi_v + \gamma \phi_{uv}) \\ \frac{-1}{\phi(X,t)}(\beta\phi_u +  \gamma\phi_{uv}) & 1 - \frac{1}{\phi(X,t)}(\beta\phi_{uv} +  \gamma\phi_{v})}
\]
Thus, letting $\psi = 1/\phi(X,t)$ for convenience, we obtain that
\[\begin{split}
\frac{\phi(X^{\beta,\gamma},t)}{ \phi(X,t)}  &=  \det \pmat{1 -\psi (\beta\phi_{uv} + \gamma \phi_{u}) & -\psi(\beta\phi_v + \gamma  \phi_{uv}) \\ -\psi(\beta\phi_u +  \gamma\phi_{uv}) & 1-\psi(\beta\phi_{uv} +  \gamma\phi_{v})} \\
&= \left(1 -\psi (\beta\phi_{uv} + \gamma \phi_{u}) \right)\left(1-\psi(\beta\phi_{uv} +  \gamma\phi_{v})\right)   \\
&\quad- \left(\psi(\beta\phi_v + \gamma  \phi_{uv}) \right) \left( \psi(\beta\phi_u +  \gamma\phi_{uv})  \right)  \\
&=  1 -2 \psi \beta\phi_{uv} -\psi \gamma \phi_{v} -\psi \gamma \phi_{u} +\psi^2 (\beta^2\phi_{uv}^2 +  \beta\gamma\phi_{v}\phi_{uv} + \beta\gamma \phi_{uv}\phi_{u} + \gamma^2 \phi_u \phi_v )  \\
&\quad-\psi^2 (\beta^2\phi_v\phi_u + \beta\gamma\phi_v \phi_{uv} + \beta\gamma \phi_{uv}\phi_u + \gamma^2 \phi^2_{uv})\\
&=  1 -2 \psi \beta\phi_{uv} -\psi \gamma \phi_{v} -\psi \gamma \phi_{u} +\psi^2 \beta^2\phi_{uv}^2+ \psi^2 \gamma^2 \phi_u \phi_v  \\
&\quad-\psi^2 \beta^2\phi_v\phi_u - \psi^2 \gamma^2 \phi^2_{uv}.\\
\end{split}
\]
We now apply Lemma \ref{lem:wab} to obtain that 
\[\begin{split}
\frac{\phi(X^{\beta,\gamma},t)}{ \phi(X,t)}
&=  1 -2 \psi \beta\phi_{uv} -\psi \gamma \phi_{v} -\psi \gamma \phi_{u} +\psi^2 \beta^2 (\phi_u\phi_v - \phi(X,t) \phi(X\diff \{u,v\},t)) \\
&\quad + \psi^2 \gamma^2 \phi_u \phi_v  
-\psi^2 \beta^2\phi_v\phi_u - \psi^2 \gamma^2  (\phi_u\phi_v - \phi(X,t) \phi(X\diff \{u,v\},t)) \\
&=  1 -2 \psi \beta\phi_{uv} -\psi \gamma \phi_{v} -\psi \gamma \phi_{u} - \psi^2 \beta^2  \phi(X,t) \phi(X\diff \{u,v\},t) \\
&\quad  + \psi^2 \gamma^2  \phi(X,t) \phi(X\diff \{u,v\},t) \\
&= \frac{1}{\phi(X,t)}\left(\phi(X,t) -2 \beta\phi_{uv} - \gamma (\phi_{v}+\phi_{u})
 +  (\gamma^2 -\beta^2) \phi(X\diff \{u,v\},t)  \right).\\
\end{split}
\]
We have obtained the following:
\begin{equation}\label{eq:phixbeta}
\phi(X^{\beta,\gamma},t) 
= \phi(X,t) -2 \beta\phi_{uv} - \gamma (\phi_{v}+\phi_{u})
 +  (\gamma^2 -\beta^2) \phi(X\diff \{u,v\},t)  .
\end{equation} 
If $u$ and $v$ are cospectral, then $\phi_u = \phi_v$ and 
\[\begin{split}
\frac{\phi(X^{\beta,\gamma},t)}{ \phi(X,t)}
&=  1 -2 \psi \beta\phi_{uv} -2\psi \gamma \phi_{u}  +\psi^2 \beta^2\phi_{uv}^2+ \psi^2 \gamma^2 \phi_u^2 
-\psi^2 \beta^2\phi_u^2 - \psi^2 \gamma^2 \phi^2_{uv}\\
&=  1 -2 \psi \beta\phi_{uv} -2\psi \gamma \phi_{u}  +\psi^2 (\beta^2-\gamma^2)\phi_{uv}^2 - \psi^2 (\beta^2-\gamma^2) \phi_u^2  \\
&=  1 -2 \psi \beta\phi_{uv} -2\psi \gamma \phi_{u}  +\psi^2 (\beta^2-\gamma^2)(\phi_{uv}^2 -  \phi_u^2 ) \\
&= \left( 1 -\psi (\beta-\gamma)(\phi_{uv}-\phi_{u})\right) \left(1 - \psi (\beta+\gamma)(\phi_{uv} + \phi_{u}) \right)
\end{split}
\]
as claimed. \qed

A graph $X$ with adjacency matrix $A$ is said to be 1-walk-regular if for all $\ell \in N$, there exist constants $a_{\ell}$ and $b_{\ell}$ such that
\begin{enumerate}[(i)]
\item $A^{\ell}\circ I = a_{\ell}I$; and
\item $A^{\ell}\circ A = b_{\ell}A$.
\end{enumerate}
Strongly regular graphs are examples of 1-walk-regular graphs. If $X$ is a 1-walk-regular graph with edges $uv$ and $xy$, then we see from the definition that 
\[W_{u,v}(X,t) = W_{x,y}(X,t)
\]
and 
\[
W_{w,w}(X,t) = W_{u,u}(X,t)
\]
for any two vertices $u,w$. We obtain the following as a corollary about 1-walk-regular graphs. 

\begin{corollary}\label{cor:1wlkreg}
If $X$ 1-walk-regular graph with edges $uv$ and $xy$, then 
\[\phi\left(X^{\beta,\gamma}_{u,v},t\right) = \phi\left(X^{\beta,\gamma}_{x,y},t\right) .
\]
In particular, when $\beta = -1$ and $\gamma =0$, we have that $X$ with edge $uv$ deleted is cospectral to $X$ with edge $xy$ deleted.
\end{corollary}

\proof It is clear that  $\phi_{v} = \phi_{u}= \phi_{x} = \phi_{y}$ and we see from the above that $\phi_{uv} = \phi(x,y)$. Thus we have that $\phi(X\diff \{u,v\},t) = \phi(X\diff \{x,y\},t) $ from Lemma \ref{lem:wab}. The statement follows from Theorem \ref{thm:beta}.\qed 

Note that we can obtain statements about the spectrum of a graph after deleting or adding an edge between $u$ and $v$ by substituting $\gamma =0$ and $\beta = \pm 1$. 

If $u$ and $v$ are cospectral in $X$, we would like to determine when they are strongly cospectral in $X^{\beta,\gamma}$. We assume, from this point, that $\phi_u = \phi_v$.  Let $U(+)$ and $U(-)$ respectively denote the $A$-modules 
generated by $e_u+e_v$ and $e_u-e_v$ and let $U(0)$ be the orthogonal complement of the 
sum $U(+)+U(-)$. Because $u$ and $v$ are cospectral, the subspaces $U(+)$ and $U(-)$
are orthogonal and $A$-invariant and thus
\[
	\re^{|V(X)|} = U(+) \oplus U(-) \oplus U(0).
\]
is an $A$-invariant decomposition of $\re^{|V(X)|}$. We can prove a stronger statement below. For cospectral vertices $u$ and $v$, we say that an $A$-invariant subspace $W$ of $\re^{|V(X)|}$ is \textsl{balanced} (resp. \textsl{skew}) with respect to $u$ and $v$ if every eigenvector $\Zw$ of $A$ such that  $\Zw \in W$ has the property that $\Zw_u = \Zw_v$ (resp.  $\Zw_u = -\Zw_v$).

\begin{lemma}\label{lem:uplusmin0}
Let $M$ be any matrix whose column space is spanned by $e_u$ and $e_v$. The subspaces $U(+)$, $U(-)$ and $U(0)$ are invariant under the algebra $\langle A,M\rangle$.
\end{lemma} 

\proof
Since the column space of $M$ is spanned by $e_u$ and $e_v$, 
it is orthogonal to $U(0)$ and therefore $MU(0)=0$. Thus
$U(0)$ is invariant under $\langle A,M\rangle$.

Since each vector in $U(+)$ is balanced relative to $u$ and $v$, each vector
in $U(+)$ is mapped by $M$ to a scalar multiple of $e_u+e_v$, and therefore $U(+)$ 
is $M$-invariant. Similarly if $z\in U(-)$, then $z$ is skew and $Mz$ is a scalar
multiple of $e_u-e_v$. So our decomposition
\[
    U(+) \oplus U(-) \oplus U(0)
\]
is invariant under both $A$ and $M$, and hence is invariant under the algebra $\langle A,M\rangle$
generated by $A$ and $M$. \qed

In particular, Lemma \ref{lem:uplusmin0} holds for $M =H$. As the orthogonal decomposition
\[
	\re^{|V(X)|} = U(+) \oplus U(-) \oplus U(0)
\]
is $A$-invariant, there is a basis of eigenvectors for $A$ such that each vector in the basis lies in one of the three summands in this
decomposition. We can decide which it is by its eigenvalue. 

\begin{lemma}\label{lem:balorskew}
	Suppose $z$ is an eigenvector for $A(X^{\beta,\gamma})$ with eigenvalue $\lambda$. Then,
	\begin{enumerate}[(a)]
		\item 
		If $z\in U(0)$, then $\lambda$ is an eigenvalue of $X$ and $X\diff\{u,v\}$.
		\item
		If $z\in U(+)$ and $\beta \neq - \gamma$, then $\lambda$ is a zero of the rational function
		\[
			1 -\frac{(\beta+\gamma)(\phi_{u,v}(X,t) + \phi(X\diff u,t))}{\phi(X,t)}  = 0.
		\]
		\item
		If $z\in U(-)$ and $\beta \neq \gamma$, then $\lambda$ is a zero of the rational function
		\[
			1 -\frac{(\beta-\gamma)(\phi_{u,v}(X,t) - \phi(X\diff u,t))}{\phi(X,t)}  = 0.
		\]		
	\end{enumerate}
\end{lemma}

\proof
Suppose $(A+H)z = \lambda z$. If $z\in U(0)$ then $Hz=0$ and therefore
$\lambda$ is an eigenvalue of $A$ (and of $A(X\diff\{u,v\}$). 

If $z\notin U(0)$, we may assume that 
\[
    z \in U(0)^\perp = U(+) + U(-).
\]
and hence that $z$ lies in $U(+)$ or $U(-)$. We suppose it is
in $U(+)$, and there is no loss in assuming additionally that $e_u^Tz=e_v^Tz=1$. Observe that $Hz = (\beta + \gamma)(e_u+e_v)$.
If $(A+ H)z = \lambda z$, then
\[
    z = (\lambda I-A)^{-1} H z =   (\lambda I-A)^{-1} (\beta + \gamma)(e_u+e_v),
\]
provided $\beta \neq - \gamma$. 
If we apply $H$ to both sides of this we get
\[\begin{split}
   H z &=  H(\lambda I-A)^{-1} (\beta + \gamma)(e_u+e_v) \\
    (\beta + \gamma)(e_u+e_v) &=  H(\lambda I-A)^{-1} (\beta + \gamma)(e_u+e_v) \\
    (\beta + \gamma)(I- H(\lambda I-A)^{-1})(e_u+e_v) &= 0,
    \end{split}
\]
which implies that
\[ \begin{split}
\left(	\pmat{1&0\\0&1} -\frac{1}{\phi(X,t)} \pmat{\beta\phi_{uv} + \gamma \phi_{u} & \beta\phi_u + \gamma  \phi_{uv} \\ \beta\phi_u +  \gamma\phi_{uv} & \beta \phi_{uv} +  \gamma\phi_{u}} \right) \pmat{1\\1}
		&= \pmat{0\\0} \\
\pmat{1\\1} -\frac{1}{\phi(X,t)} \pmat{(\beta+\gamma)\phi_{uv} + (\beta+\gamma) \phi_{u}  \\  (\beta + \gamma) \phi_{uv} + (\beta+\gamma)\phi_u  }  
		&= \pmat{0\\0} .
		\end{split}
\]
Thus, we obtain that $\lambda$ must be a root of the rational function
\[
1 -\frac{(\beta+\gamma)\phi_{uv} + (\beta+\gamma)\phi_{u}}{\phi(X,t)}  = 0.
\]
If $z\in U(-)$ and $e_u^Tz= - e_v^Tz=1$. Observe that $Hz = (\gamma - \beta)(e_u-e_v)$.
A similar argument shows that in this case, we have, provided $\beta \neq \gamma$, that 
\[\left(	\pmat{1&0\\0&1} -\frac{1}{\phi(X,t)} \pmat{\beta\phi_{uv} + \gamma \phi_{u} & \beta\phi_u + \gamma  \phi_{uv} \\ \beta\phi_u +  \gamma\phi_{uv} & \beta \phi_{uv} +  \gamma\phi_{u}} \right) \pmat{1\\-1}
		= \pmat{0\\0} \\
\] and thus  $\lambda$ must be a zero of
\[
1 -\frac{(\beta-\gamma)\phi_{uv} - (\beta -\gamma)\phi_{u}}{\phi(X,t)}  = 0,
\]
and the lemma follows.\qed

\begin{theorem} \label{thm:perturb-str-cospec}
    Suppose $u$ and $v$ are cospectral vertices in $X$.
\begin{enumerate}[(a)]
\item If $\beta\neq\pm\gamma$, then $u$ and $v$ 
	are strongly cospectral in $X^{\beta,\gamma}$ if and only if the rational functions
    \[
        1-\frac{(\beta+\gamma)(\phi_{u,v}(X,t)+\phi(X\diff u,t))}{\phi(X,t)},
        \,\,
            1-\frac{(\beta-\gamma)(\phi_{u,v}(X,t)-\phi(X\diff u,t))}{\phi(X,t)}
    \]
    have no common zeros.
\item If $\beta = \gamma$ but $\beta\neq-\gamma$, then $u$ and $v$ are strongly cospectral in $X^{\beta,\gamma}$ if and only if 
\[
        1-\frac{(\beta+\gamma)(\phi_{u,v}(X,t)+\phi(X\diff u,t))}{\phi(X,t)},
        \,\,
            \phi(X,t)
    \]
have no common zeros.
\item If $\beta=-\gamma$ but $\beta\neq\gamma$, then $u$ and $v$ are strongly cospectral in $X^{\beta,\gamma}$ if and only if
\[
         1-\frac{(\beta-\gamma)(\phi_{u,v}(X,t)-\phi(X\diff u,t))}{\phi(X,t)},
        \,\,
           \phi(X,t)
    \]
have no common zeros.
\end{enumerate}
\end{theorem}  

\proof
Every eigenvector of $A + H$ in $U(+)$ is balanced relative to $u$ and $v$, 
while those in $U(-)$ are skew relative to $u$ and $v$.  If $\beta=\gamma$, then it is clear that an eigenvector of $A+H$ from $U(-)$ is an eigenvector if $A$, and likewise, if $\beta=-\gamma$, then an eigenvector of $A+H$ from $U(+)$ is an eigenvector of $A$.   So in each case, if these rational functions have no common zeros, then Theorem \ref{thm:beta}, and Lemma~\ref{lem:balorskew} give us that each eigenspace is either balanced or skew. This in turn implies that $E_re_u=\pm E_re_v$ for each idempotent $E_r$, which implies the theorem.\qed

\section{Perturbation in strongly regular graphs}\label{sec:srg}

For this section, we will consider $X$, a strongly regular graph of valency $k$ on $n$ vertices with eigenvalues $k$, $\theta$ and $\tau$, with multiplicities $1, m_\theta,$ and $m_\tau$ respectively. The idempotent matrices in the spectral decomposition of $A := A(X)$ can be written as linear combinations of $I$, $A$, and the adjacency matrix of the complement of $X$, which we will denote $\comp{A}$. The spectral decomposition of $A$ is 
\[
A = kE_k + \theta E_{\theta} + \tau E_{\tau},
\]
where 
\[
\begin{split}
E_k &= \frac{1}{n} J =  \frac{1}{n}(I + A + \comp{A}), \\
E_{\theta} &= \frac{m_{\theta}}{n}\left(I +\frac{\theta}{k} A - \frac{\theta + 1}{n-k-1} \comp{A}\right), \text{ and,} \\
E_{\tau} &= \frac{m_{\tau}}{n}\left(I +\frac{\tau}{k} A - \frac{\tau + 1}{n-k-1} \comp{A}\right).
\end{split}\]

\begin{theorem}\label{thm:srgperturb}
If $u$ an $v$ are adjacent vertices of $X$, then the eigenvalues of the perturbation of $X$ at vertices $u,v$ with weights $\beta,\gamma$ are as follows: 
\begin{enumerate}[(a)]
\item 
    $\theta$ with multiplicity $m_{\theta}-2$ and $\tau$ with multiplicity $m_{\tau}-2$.
    \item
    The five solutions of the following two rational equations 
    \begin{align}
    \frac{2k}{t-k} +\frac{(k+\theta) m_{\theta}}{t-\theta} +  \frac{(k+\tau) m_{\tau}}{t-\tau}
        &= \frac{nk}{\beta+\gamma} \label{eq:srgab1}\\
     \frac{(k-\theta) m_{\theta}}{t-\theta} +  \frac{(k-\tau) m_{\tau}}{t-\tau} &= \frac{nk}{\beta-\gamma}. \label{eq:srgab2}
\end{align}
\end{enumerate}
If $u$ an $v$ are non-adjacent vertices of $X$, then the eigenvalues of the perturbation of $X$ at vertices $u,v$ with weights $\beta,\gamma$ are as follows:
\begin{enumerate}[(a)]
  \setcounter{enumi}{2}
\item 
    $\theta$ with multiplicity $m_{\theta}-2$ and $\tau$ with multiplicity $m_{\tau}-2$.
    \item
    The five solutions of the following two rational equations 
    \begin{align} 
    \frac{2(n\!-\!k\!-\!1)}{t-k} \!+\!\frac{(n\!-\!k\!-\!2\!-\!\theta) m_{\theta}}{t-\theta} \!+\!  \frac{(n\!-\!k\!-\!2\!-\!\tau) m_{\tau}}{t-\tau}
        &= \frac{n(n\!-\!k\!-\!1)}{\beta+\gamma} \label{eq:srgab3}\\
     \frac{(n-k-\theta ) m_{\theta}}{t-\theta} +  \frac{(n-k- \tau) m_{\tau}}{t-\tau} &= \frac{n(n\!-\!k\!-\!1)}{\beta-\gamma}. \label{eq:srgab4}
\end{align}
\end{enumerate}
\end{theorem}

\proof If $u$ and $v$ are adjacent vertices of $X$, we have 
\[
\frac{\phi_{uv}(X,t)}{\phi(X,t)} = \frac{1}{n} \left( \frac{1}{t-k} +  \frac{m_{\theta}}{t-\theta}\frac{\theta}{k} +  \frac{m_{\tau}}{t-\tau} \frac{\tau}{k} \right)
 = \frac{1}{nk} \left( \frac{k}{t-k} +  \frac{\theta m_{\theta}}{t-\theta} +  \frac{\tau m_{\tau}}{t-\tau}  \right)
\]
from (\ref{eq:phiabEr}) and 
\[
\frac{\phi(X\diff  u,t)}{\phi(X,t)} = \frac{1}{n} \left( \frac{1}{t-k} +  \frac{m_{\theta}}{t-\theta} +  \frac{m_{\tau}}{t-\tau}\right)
\]
from (\ref{eq:phiaEr}). Reprising our abbreviations from the previous section, we see that 
\[ \begin{split}
\frac{\phi_u-\phi_{uv}}{\phi(X,t)} &= \frac{1}{kn} \left( \frac{k m_{\theta}}{t-\theta} -   \frac{\theta m_{\theta}}{t-\theta} +  \frac{k m_{\tau}}{t-\tau} - \frac{\tau m_{\tau}}{t-\tau} \right) \\
&= \frac{1}{kn} \left( \frac{(k-\theta) m_{\theta}}{t-\theta} +  \frac{(k-\tau) m_{\tau}}{t-\tau}\right)
\end{split}
\]
and 
\[
\frac{\phi_u+\phi_{uv}}{\phi(X,t)} 
= \frac{1}{kn} \left(  \frac{2k}{t-k} +\frac{(k+\theta) m_{\theta}}{t-\theta} +  \frac{(k+\tau) m_{\tau}}{t-\tau}\right).
\]

Since every pair of vertices in a strongly regular graph is cospectral, we have from Theorem \ref{thm:beta} that 
\[ \begin{split}
\frac{\phi(X^{\beta,\gamma},t)}{ \phi(X,t)} &= \left( 1 -\frac{(\beta-\gamma)(\phi_{uv}-\phi_{u})}{\phi(X,t)}\right) \left(1 - \frac{(\beta+\gamma)(\phi_{uv} + \phi_{u}) }{\phi(X,t)} \right) \\
&= \left(1- \frac{\beta - \gamma}{kn} \left( \frac{(k-\theta) m_{\theta}}{t-\theta} +  \frac{(k-\tau) m_{\tau}}{t-\tau}\right) \right) \\
& \quad * \left(1- \frac{\beta+\gamma}{kn} \left(  \frac{2k}{t-k} +\frac{(k+\theta) m_{\theta}}{t-\theta} +  \frac{(k+\tau) m_{\tau}}{t-\tau}\right) \right).
\end{split}\]

This leads to equations \eqref{eq:srgab1} and \eqref{eq:srgab2} and we obtain that 
the eigenvalues of $X^{\beta,\gamma}$ are as stated in the lemma. 

If $u$ and $v$ are not adjacent, a similar argument will give equations  \eqref{eq:srgab3} and \eqref{eq:srgab4} and the statement follows. \qed 

We obtain the following immediate corollary.

\begin{corollary}
For any strongly regular graph $X$ and any edge $e$ of $X$, the eigenvalues of $X$ with edge $e$ deleted does not depend on the choice of $e$. Similarly, let $u,v$ be non-adjacent vertices of $X$. The eigenvalues of the graph $X$ with an edge added joining $u$ and $v$ does not depend on the choice of $u$ and $v$. 
\end{corollary}

Observe that if the five `new' zeros are distinct from one another, then $u$ and $v$
are strongly cospectral in the perturbed graph. If the new zeros are not distinct,
they overlap at a zero of $\phi(X\diff u)$, equivalently at a zero of 
$\phi'(X,t)$.  

We give a lemma on strongly regular graphs which we will use to simplify these polynomials further.
\begin{lemma}\label{lem:srg_identites}
Let $X$ be a strongly regular graph with parameter $(n,k,\theta,\tau)$, and non-trivial eigenvalues $\theta,\tau$ with multiplicities $m_\theta,m_\tau$ respectively.  The following identities hold:
\begin{align}
    m_\theta+m_\tau+1 &= n\label{mult}\\
    \theta m_\theta + \tau m_\tau + k &= 0\label{trace}\\
    \theta^2m_\theta + \tau^2m_\tau + k^2 &= nk\label{trace2}\\
    \theta m_\tau + \tau m_\theta &= (n-1)(\theta+\tau)+k\label{mix}\\
    \theta\tau m_\theta + \theta\tau m_\tau &= k^2-nk-k(\theta+\tau)\label{mix2}\\
    (\theta\tau +k)n &= k^2-k(\theta+\tau)+\theta\tau.\label{srg_identity}
\end{align}
\end{lemma}
\proof
Identities (\ref{mult}), (\ref{trace}), and (\ref{trace2}) follow form considering the trace of $A^r$ for $r=0,1,2$ respectively.  Then (\ref{mix}) and (\ref{mix2}) follow by manipulation of these.

Identity (\ref{srg_identity}) follows from the well-known (see \cite{GodsilRoyle}) facts that the parameters of a strongly regular graph satisfy
\[
(n-k-1)c = k(k-a-1)
\]
and that $\theta$ and $\tau$ are roots of the quadratic $t^2-(a-c)t-(k-c).$
\qed

\begin{corollary}\label{lem:cubic-quadratic-equations}
Let $X$ be a strongly regular graph with eigenvalues $k$, $\theta$ with multiplicity $m_\theta$, and $\tau$ with multiplicity $m_\tau$ as above.  Then the eigenvalues of $X^{\beta,\gamma}$ are $\theta$ with multiplicity $m_\theta-2$, $\tau$ with multiplicity $m_\tau-2$, and the roots of the equations given as follows:
\begin{align} \label{eq:cubic-adj}
    (t\!-\!k)(t\!-\!\theta)(t\!-\!\tau) &= (\beta\!+\!\gamma)\left(t^2-(k\!+\!\theta\!+\!\tau\!-\!1)t + k\theta\!+\!k\tau\!+\!k\!+\!2\theta\tau\right) \\ \label{eq:quadratic-adj}
    (t-\theta)(t-\tau) &= (\beta-\gamma)\left(t-(\theta+\tau+1)\right)
\end{align}
when $u$ and $v$ are adjacent; and
\begin{align} \label{eq:cubic-nonadj}
    (t\!-\!k)(t\!-\!\theta)(t\!-\!\tau) &= (\beta\!+\!\gamma)\left(t^2-(k\!+\!\theta\!+\!\tau)t + k\theta\!+\!k\tau\!+\!2k\!+\!2\theta\tau\right)\\ \label{eq:quadratic-nonadj}
    (t-\theta)(t-\tau) &= (\beta-\gamma)\left(t-(\theta+\tau)\right)
\end{align}
when $u$ and $v$ are not adjacent.
\end{corollary}

\proof
First let us suppose that $u$ and $v$ are adjacent.  We have that the new eigenvalues for $X^{\beta,\gamma}$ are 
\begin{align*}
    \frac{2k}{t-k}+\frac{(k+\theta)m_\theta}{t-\theta}+\frac{(k+\tau)m_\tau}{t-\tau} &= \frac{nk}{\beta+\gamma}\\
    \frac{(k-\theta)m_\theta}{t-\theta} + \frac{(k-\tau)m_\tau}{t-\tau} &= \frac{nk}{\beta-\gamma}.
\end{align*}

We will make use of the identities from Lemma \ref{lem:srg_identites}. Let us first work with the quadratic.  In the rational expression above, if we clear denominators we obtain
\[
nk(t-\theta)(t-\tau) = (\beta-\gamma)\left((k-\theta)m_\theta(t-\tau)+(k-\tau)m_\tau(t-\theta)\right).
\]
Multiplying the right side out and applying identities (\ref{mult}), (\ref{trace}), (\ref{mix}), and (\ref{mix2}) yields
\[
(t-\theta)(t-\tau) = (\beta-\gamma)\left(t-(\theta+\tau+1)\right)
\]

For the cubic, clearing denominators 
and multiplying out the right side, the coefficient of $t^2$ and $t$ can be obtained in much the same way that we did the quadratic above. 

The constant term on the right side becomes 
\[
(\beta+\gamma)\left(2k\theta\tau+(k+\theta)m_\theta \tau k + (k+\tau)m_\tau \theta k\right).
\]
Applying (\ref{mix}) and (\ref{mix2}), this becomes
\[
(\beta+\gamma)\left(2k\theta\tau + k^2((n-1)(\theta+\tau)+k)+k(k^2-nk-k(\theta+\tau))\right).
\]
Now we will apply the identity (\ref{srg_identity}) of Lemma \ref{lem:srg_identites}.
The above becomes
\[
(\beta+\gamma)k\left(2(\theta\tau +k)n+kn(\theta+\tau)-nk\right)
\]
which becomes
\[
(\beta+\gamma)nk(2\theta\tau + k + k\theta+k\tau)
\]
and the result follows.

The case where $u$ and $v$ are non-adjacent is done similarly. 
\qed

\section{Perfect state transfer}\label{sec:pst}
In this section we prove that under certain conditions a perturbation of $X$ at vertices $u,v$ with weights $\beta,\gamma$ of a strongly regular graph exhibits perfect state transfer. We will continue with the notation from the previous section. 

Recall that a necessary condition for pst between $u$ and $v$ to occur is that $u$ and $v$ must be strongly cospectral; that is, $E_re_u = \pm E_re_v$ for all spectral idempotents $E_r$. When $u$ and $v$ are strongly cospectral, the eigenvalues whose eigenvectors are supported on $u$ and $v$ can be naturally grouped into two sets:  we will let $\mu_1,...,\mu_\ell$ denote eigenvalues with $E_{\mu_i}e_u = E_{\mu_i}e_v$ and $\lambda_1,...,\lambda_m$ denote the eigenvalues for which $E_{\lambda_i}e_u = -E_{\lambda_i}e_v$.

We will make use of the following lemma to study perfect state transfer on a graph by investigating its eigenvalues and eigenvectors (see \cite{CoutinhoGodsilGuoVanhove2,KemptonLippnerYauPath}).  
For an integer $m$ we will let $|m|_2$ denote the highest power such that $2^{|m|_2}$ divides $m$.

\begin{lemma}\cite{CoutinhoGodsilGuoVanhove2}\label{lem:pst_eig}
Perfect state transfer between vertices $u$ and $v$ of $X$ occurs at time $t$ if and only if the following two conditions are satisfied:
\begin{enumerate}[(i)]
\item the vertices $u$ and $v$ are strongly cospectral; and
\item for each $i$,
\[
\frac{\lambda_1 - \mu_i}{\lambda_1 - \lambda_2} = \frac{p_i}{q_i} \in \Q
\]
where $p_i$ is an odd integer and $q_i$ is an even integer such that $|q_i|_2=r$ for all $i$ for some constant $r$.
\end{enumerate}
\end{lemma}

In this section, we will show our main result on perfect state transfer in perturbed strongly regular graphs, which is summarized in Theorem \ref{thm:pst}. 

We will proceed by first introducing some notation, which will be used throughout this section. Let  $u,v$ be a pair of non-adjacent vertices. Let us assume $\beta = -\gamma$, so the roots of the cubic equation \eqref{eq:cubic-nonadj}  are simply $\theta, \tau, k$.  Let us write $\beta - \gamma = p/q$ where $p,q$ are coprime integers. 

Let us denote by $\lambda_{1,2}$ the roots of \eqref{eq:quadratic-nonadj}. Then 
\[ \lambda_{1,2} = \frac{p/q + \theta+\tau \pm \sqrt{(p/q + \theta + \tau)^2 - 4 (p/q) (\theta+\tau) - 4\theta \tau}}{2}\]
Let us introduce the notation $\alpha/q$ for the discriminant, that is 
\[ \frac{\alpha}{q} = \sqrt{(p/q + \theta + \tau)^2 - 4 (p/q) (\theta+\tau) - 4\theta \tau}.\] That is,
\[ \alpha^2 = (p - q(\theta+\tau))^2 - 4q^2 \theta \tau,\] or equivalently 
\[ 4q^2 \theta \tau = (p - q(\theta + \tau) - \alpha)(p - q(\theta+\tau) +\alpha).\]

Let us further write $A = p - q(\theta + \tau)-\alpha$ and $B = p - q(\theta + \tau) + \alpha$. Then 
\begin{enumerate}[(i)]
\item $2\alpha = B - A$;
\item $4q^2 \theta \tau = A \cdot B$; and
\item $\lambda_1 = \frac{B + 2q(\theta+\tau)}{2q}$.
\end{enumerate}

If $\mu$ denotes any one of $\theta, \tau, k$, then we get 
\begin{equation}\label{eq:odd-even} \frac{\lambda_1 - \mu}{\lambda_1 - \lambda_2} = \frac{B + 2q(\theta+\tau - \mu)}{B - A}. \end{equation}

It is this quotient that has to be odd over even for each value of $\mu$, moreover the even denominator needs to have the same power of two in order to have perfect state transfer. 

 Let us note some immediate consequences of this odd-even condition:
\begin{enumerate}[(i)]
\item Since $2\alpha = A-B$, we get $A \equiv B \pmod{2}$. 
\item Since 4 divides $A \cdot B$, both $A$ and $B$ have to be even.
\item The numerator of \eqref{eq:odd-even} is also even and so the denominator must be divisible by 4. Thus we have that $\alpha$ is also even.
\end{enumerate}

Now suppose $q$ is even. Since $p$ and $q$ are coprime, then $p$ must be odd. However, since $A = p-q(\theta+\tau) - \alpha$ is also even, and $q$ is even, then $p+\alpha$ must be even, which would mean that $\alpha$ is odd, contradicting (iii) above. 

So $q$ must be odd. We know $A,B$ must both be even, thus we write $A = 2A', B = 2B'$. Then $B' - A' = \alpha$ and $q^2 \theta \tau = A' B'$. Furthermore we get the new conditions that 
\begin{equation}\label{eq:odd-even2} \frac{\lambda_1 - \mu}{\lambda_1 - \lambda_2} = \frac{B' + q(\theta+\tau - \mu)}{B' - A'} \end{equation}
must be odd over even with the even denominator having the same power of 2 for each of $\mu = \theta, \tau, k$. This in particular implies that $\theta \equiv \tau \equiv k \pmod{2}$ is required. 

Suppose first that all three are odd. Then $A' \cdot B'$ is odd hence both $A'$ and $B'$ are odd. But then the numerator in \eqref{eq:odd-even2} is still even, so we need the denominator to be divisible by 4. The denominator is $B'- A'$, so we need $A' \equiv B' \pmod{4}$, which is equivalent to requiring $A' \cdot B' = q^2 \theta \tau \equiv 1 \pmod{4}$. Hence we need $\theta \tau \equiv 1 \pmod{4}$, which is again equivalent to asking $\theta \equiv \tau \pmod{4}$. However we also need that the numerator in \eqref{eq:odd-even2} is divisible by the same power of 2 for all three choices of $\mu$, which also implies that $\theta \equiv k \pmod{4}$ is a requirement. 
 
In the lemma below, note that we may choose $q$ such that $q \equiv \theta \mod 4$ and we may also choose any factorization of $q = q_1q_2$, which determines the choice of $p$. The values of $p/q$, which $p,q$ is a valid choice, is dense in $\Q$.

\begin{lemma}\label{lem:claimodd}
If $\theta \equiv \tau \equiv k \pmod{4}$ are odd, then for any $q \equiv \theta \mod 4$, and $p=q^2+ \theta \tau+q(\theta+\tau)$, setting $\beta-\gamma = p/q$ implies the three quotients of \eqref{eq:odd-even} satisfy the parity conditions. In fact, the set of values $p/q$ for which the parity conditions are satisfied with $\beta-\gamma = p/q$ is dense in $\Q$. 
\end{lemma}

\proof
Let's choose $B' = \theta \tau$ and $A' = q^2$. This means $p = q^2+ \theta \tau+q(\theta+\tau)$. Since $\theta \equiv \tau \pmod{4}$, we have $B' \equiv A' \equiv 1 \pmod{4}$ and thus $B' - A' \equiv 0 \pmod{4}$. On the other hand the numerator of \eqref{eq:odd-even2} is $B' + q(\theta+\tau-\mu) \equiv 1 + q \theta \pmod{4}$. We need this to be congruent to $2 \pmod{4}$, which can be always achieved by choosing $q \equiv \theta \pmod{4}$.

To see the density claim, let's write $q= q_1 q_2$ and choose $A' = q_1^2$ and $B' = q_2^2 \theta \tau$. Then $p = q_1^2+q_2^2\theta \tau + q_1 q_2 (\theta+\tau)$. Then, given these choices,
\[ \frac{p}{q} = \frac{q_1^2+q_1 q_2(\theta+\tau)+q_2^2\theta \tau}{q_1 q_2} = \theta + \tau + \frac{q_2}{q_1}\theta \tau + \frac{q_1}{q_2}.\]

Any choice of $q_1, q_2$ is valid as long as $q_1 q_2 \equiv \theta \pmod{4}$, so this still allows a dense set of rational $q_1/q_2$ values. Thus, since $\theta \tau < 0$ always, the set of valid values of $p/q$ is also dense in $\rats$.
\qed

Note that the assumptions of Lemma \ref{lem:claimodd} imply that $u$ and $v$ are strongly cospectral; this follows from Theorem \ref{thm:perturb-str-cospec} and from the fact that the $U(+)$ eigenspace has roots $\lambda_1, \lambda_2$, and by the parity condition in Lemma \ref{lem:pst_eig}, these roots will automatically be disjoint from any possible root ($k, \theta, \tau$) of the $U(-)$ eigenspace. 

\begin{lemma}\label{thm:pstsrg-nonadj}
There is perfect state transfer between a pair of non-adjacent nodes in an SRG with a perturbation of $X$ at vertices $u,v$ with weights $\beta,\gamma$ of the form $\beta + \gamma = 0$ such that $\beta$ is rational, if and only if there are odd integers $k', \theta', \tau'$ congruent to each other modulo 4 and another non-negative integer $t$ such that $k = 2^t k', \theta = 2^t \theta', \tau = 2^t \tau'$. If this is the case, then there is perfect state transfer from $u$ to $v$ in the perturbation of $X$ at vertices $u,v$ with weights $\beta,\gamma$ where $\beta = \theta'+\tau'$.

In particular this always holds when $k \equiv \theta \equiv \tau \pmod{4}$ are odd numbers.
\end{lemma}

\proof First let us consider when $\theta, \tau, k$ are all even. In this case it is easy to see that both $A', B'$ need to be even, but then we can divide each of them by 2 and get a new set of integer values that again have to satisfy the same parity conditions as in \eqref{eq:odd-even2}. This can be repeated until one of $\theta, \tau$, and $k$ becomes odd, at which moment all three must become odd, and the the lemma follows from Lemma \ref{lem:claimodd}. The particular choice of $\beta$ comes from choosing $q=\theta'$ in Lemma \ref{lem:claimodd}. \qed 

\begin{lemma} \label{thm:pstsrg-adj}
There is perfect state transfer between a pair of adjacent nodes in an SRG with a single-edge perturbation of the form $\beta + \gamma = 0$ such that $\beta$ is rational, if and only if there are odd integers $k', \theta', \tau'$ congruent to each other modulo 4 and another non-negative integer $t$ such that $k+1 = 2^t k', \theta+1 = 2^t \theta', \tau+1 = 2^t \tau'$. If this is the case, then there is perfect state transfer from $u$ to $v$ in the perturbation of $X$ at vertices $u,v$ with weights $\beta,\gamma$ where $\beta =  \theta'+\tau'$.

In particular this always holds when $k \equiv \theta \equiv \tau \pmod{4}$ are even numbers.
\end{lemma}

\proof We consider the case where $k \equiv \theta \equiv \tau \pmod{4}$  are even numbers. We choose a pair of adjacent nodes $u,v$, and still enforce $\beta +\gamma = 0$ and denote $\beta- \gamma = p/q$. Then the roots of the cubic polynomial \eqref{eq:cubic-adj} are still $k, \theta, \tau$. The roots of the quadratic polynomial \eqref{eq:quadratic-adj} are 
\[ \lambda_{1,2} = \frac{p/q + \theta+\tau \pm \sqrt{(p/q + \theta + \tau)^2 - 4 (p/q) (\theta+\tau+1) - 4\theta \tau}}{2}.\]
Then, again denoting the discriminant by $\alpha/q$ we get 
\[ \frac{\alpha}{q} = \sqrt{(p/q + \theta + \tau)^2 - 4 (p/q) (\theta+\tau+1) - 4\theta \tau}.\] That is,
\[ \alpha^2 = (p - q(\theta+\tau+2))^2 - 4q^2 (\theta \tau +\theta+\tau+1),\] or equivalently 
\[ 4q^2 (\theta+1)( \tau+1) = (p - q(\theta+1 + \tau+1) - \alpha)(p - q(\theta+1+\tau+1) +\alpha).\]
Denoting the two factors by $A$ and $B$ respectively, we have 
\begin{enumerate}[(i)]
\item $2\alpha = B - A$;
\item $4q^2 (\theta+1) (\tau+1) = A \cdot B$; and 
\item $\lambda_1 = \frac{B + 2q(\theta+\tau+1)}{2q}$.
\end{enumerate}
Finally we obtain that 
\begin{equation}\label{eq:odd-even-adj} \frac{\lambda_1 - \mu}{\lambda_1 - \lambda_2} = \frac{B + 2q((\theta+1)+(\tau+1) - (\mu+1))}{B - A}. \end{equation}
From here, it is clear that if we denote $\tilde{\theta} = \theta+1, \tilde{\tau} = \tau+1, \tilde{k} = k+1, \tilde{\mu}=\mu+1$, then the condition is that \eqref{eq:odd-even-adj} needs to be odd-over-even for any choice of $\tilde{\mu} \in \{ \tilde{\theta}, \tilde{\tau}, \tilde{k}\}$. Since now $A \cdot B = \tilde{\theta} \cdot \tilde{\tau}$, we are in an identical situation to the  one in the previous  analysis, after substituting the variables $\tilde{\theta}, \tilde{\tau}, \tilde{k}$ for $\theta, \tau, k$, respectively, and the result follows. \qed

We summarize our results in the following theorem. Note for any choice of adjacent or non-adjacent $u$ and $v$ a strongly regular graph satisfying the respective hypothesis of the following corollary,  there is a dense set of rational values for the quantityt $\beta - \gamma$ for which the the perturbation of $X$ at vertices $u,v$ with weights $\beta,\gamma$, where $\beta + \gamma = 0$, admits perfect state transfer from $u$ to $v$. Specific values of $\beta$ can be found depending on the parameters, as described in Lemmas \ref{thm:pstsrg-nonadj} and \ref{thm:pstsrg-adj}.

\begin{theorem} \label{thm:pst}
In a strongly regular graph, if $k \equiv \theta \equiv \tau \pmod{4}$ is odd, then for any pair of non-adjacent nodes $u,v$, the perturbation of $X$ at vertices $u,v$ with weights $\beta,\gamma$ admits perfect state transfer from $u$ to $v$ for some value of $\beta = -\gamma$. If $k \equiv \theta \equiv \tau \pmod{4}$ is even, then for any pair of adjacent nodes $u,v$, the perturbation of $X$ at vertices $u,v$ with weights $\beta,\gamma$ where admits perfect state transfer from $u$ to $v$ for some value of $\beta = -\gamma$.
\end{theorem}

\proof This follows directly from Lemmata \ref{thm:pstsrg-nonadj} and \ref{thm:pstsrg-adj}. \qed

\section{Examples}\label{sec:example}

In this section, we give strongly regular graphs $X$ where there exists vertices $u$ and $v$ such that  the perturbation of $X$ at vertices $u,v$ with weights $\beta,\gamma$ admits perfect state transfer from $u$ to $v$. 

An \textsl{orthogonal array} with parameters $k$ and $n$, denoted, $\OA(k,n)$, is a $k\times n^2$ array with entries from $[n]$ such that, for any two rows of the array, the $n^2$ ordered pairs of elements are all distinct. Given $\OA(k,n)$, we define a graph as follows: the vertices are the $n^2$ columns of the array and two vertices are adjacent if they are equal in exactly one coordinate position. The graph given by $\OA(k,n)$ is strongly regular with parameters
\[ (n^2, (n-1)k, n-2 + (k-1)(k-2), k(k-1))\]
and eigenvalues $(n-1)k$, $n-k$ and $-k$.
A \textsl{Latin square graph} is a graph given by an orthogonal array  $\OA(3,n)$. More background information can be found in \cite[Chapter 10]{GodsilRoyle}. 

\begin{theorem}\label{thm:oapst} Let $X$ be a graph given by $\OA(k,n)$. Let $u$ and $v$ be two distinct vertices in $X$. The perturbation of $X$ at vertices $u,v$ with weights $\beta,\gamma$ with $\beta + \gamma = 0$ and $\beta \in \rats$ admits perfect state transfer from $u$ to $v$  if and only if 
$|n|_2 \geq |k|_2 + 2$ when $u \nsim v$ and $|n|_2 \geq |k-1|_2 + 2$  when $u \sim v$.
\end{theorem}

\proof For an integer $x$, let \[f(x) = \frac{x}{2^{|x|_2}}.\] Suppose $u$ and $v$ are not adjacent. Theorem \ref{thm:pstsrg-nonadj} gives us that there is perfect state transfer from $u$ to $v$ in $X^{\beta,\gamma}$ if and only if 
\begin{equation}\label{eq:oa1} |(n-1)k|_2 =|n-k|_2 = |-k|_2
\end{equation}
and
\begin{equation}\label{eq:oa2}
f((n-1)k) \equiv f(n-k) \equiv f(-k) \pmod{4}.
\end{equation}
Let $t = |-k|_2$ and $f(k) = j$. Suppose that \eqref{eq:oa1} and \eqref{eq:oa2} are satisfied. Since $|n-k|_2 = |-k|_2$, we have that $2^t | n$ and we may write $n = 2^t q$ for some $q \in \ints$. Since $f(n-k) \equiv f(-k) \pmod{4}$, we see that $q -j \equiv -j \pmod{4}$ and thus $q \equiv 0 \pmod{4}$.

Conversely, it is clear that if $|n|_2 \geq t+2$, then \eqref{eq:oa1} and \eqref{eq:oa2} are satisfied.

 Now suppose $u$ and $v$ are adjacent in $X$. Theorem \ref{thm:pstsrg-adj} gives us that there is perfect state transfer from $u$ to $v$ in $X^{\beta,\gamma}$ if and only if 
\begin{equation}\label{eq:oa3} |(n-1)k+1|_2 =|n-k+1|_2 = |-k+1|_2
\end{equation}
and
\begin{equation}\label{eq:oa4}
f((n-1)k+1) \equiv f(n-k+1) \equiv f(-k+1) \pmod{4}.
\end{equation}
Let $t = |-k+1|_2$ and $f(k-1) = j$. Suppose that \eqref{eq:oa3} and \eqref{eq:oa4} are satisfied. Since $|n+(-k+1)|_2 = |-k+1|_2$, we have that $2^t | n$ and we may write $n = 2^t q$ for some $q \in \ints$. Since $f(n-k+1) \equiv f(-k+1) \pmod{4}$, we see that $q -j \equiv -j \pmod{4}$ and thus $q \equiv 0 \pmod{4}$.

Conversely, it is clear that if $|n|_2 \geq t+2$, then \eqref{eq:oa3} and \eqref{eq:oa4} are satisfied.
\qed

An elliptic affine polar graph of order $(2e,q)$, denoted $VO^-(2e,q)$, is the point graph of an elliptic affine polar space. These graphs give an infinite family of strongly   regular graphs with the following parameters:
\[ \left(q^{2e},\, (q^{e -1} - 1)(q^e +1),\, q(q^{e -2}-  1)(q^{e-1} +1)+q-2,\, q^{e-1} (q^e +1) \right) 
\]
where $q$ is a prime power. The two smaller distinct eigenvalues of $VO^-(2e,q)$ are $q^{e-1} - 1$ and  $-q^e +q^{e-1}-1$. For example, the 
Clebsch graph with parameters $(16,5,0,2)$ is isomorphic to $VO^-(4,2)$.
Similarly, the hyperbolic affine polar graph of order $(2e,q)$, denoted $VO^-(2e,q)$, is the point graph of a hyperbolic affine polar space and is a strongly regular graph with parameters as follows:
\[ \left(q^{2e},\, (q^{e -1} + 1)(q^e -1),\, q(q^{e -2} + 1)(q^{e-1} -1)+q-2,\, q^{e-1} (q^e +1) \right) 
\]
where $q$ is a prime power. The local structure of these graphs are studied in \cite{BrouwerShult90} and they also occur as Cayley graphs in \cite{CalderbankKantor86}. 

\begin{theorem}\label{thm:vo-pst} Let $X$ be an elliptic or hyperbolic affine polar graph and let $u$ and $v$ be distinct, non-adjacent vertices in $X$. The perturbation of $X$ at vertices $u,v$ with weights $\beta,\gamma$ with $\beta + \gamma = 0$ and $\beta \in \rats$ admits perfect state transfer from $u$ to $v$  if 
$q$ is even. 
\end{theorem}

\proof
Suppose $X$ is an elliptic affine polar graph. The eigenvalues of $X$ are 
\[
k =  (q^{e -1} - 1)(q^e +1), \theta = q^e- q^{e-1} - 1, \text{ and } \tau = - q^{e-1} - 1.
\]
We have the following:
\[
\begin{split}
k- \theta & = q^{2e -1}  - 2q^e + 2q^{e-1} \\
k - \tau &=  q^{2e -1}  - q^e + 2q^{e-1} \\
\theta - \tau &= q^e.\\
\end{split}
\]
If $q$ is even and $e > 1$, then we have that $k \equiv \theta \equiv \tau  \mod 4$ and each of $k,\theta,\tau$ are odd integers. It then follows from Theorem \ref{thm:pstsrg-nonadj} that $X$ admits perfect state transfer from $u$ to $v$. The case when $X$ is a hyperbolic affine polar graph follows similarly.
\qed

Theorems \ref{thm:oapst} and \ref{thm:vo-pst} gives infinite families of strongly regular graphs whose perturbations  at vertices $u,v$ with weights $\beta,\gamma$ admit $uv$-pst. The following theorem allows us to consider several important families of strongly regular graphs \cite{Neumaier}. 

\begin{theorem}[Neumaier, Theorem 5.1 and 5.2 of \cite{Neumaier}] 
Let $\tau$ be a fixed real number. All but finitely many strongly regular graphs with least eigenvalue $\tau$  are  conference graphs or given by a Steiner triple system or $\OA(3,n)$.
\end{theorem}

\begin{theorem}\label{srg:pst}
Let $\tau$ be a fixed real number. Besides Latin square graphs, there are at most finitely many strongly regular graphs $X$ with least eigenvalue $\tau$ such that $X$ contains vertices $u$ and $v$ such that the perturbation of $X$ at vertices $u,v$ with weights $\beta,\gamma$ with $\beta + \gamma = 0$ and $\beta \in \rats$ admits perfect state transfer from $u$ to $v$.
\end{theorem}

\proof
Strongly regular graphs given by Steiner triple systems with parameter $v$ have eigenvalues $k=3(v-3)/2$, $\theta = (v-9)/2$, and $\tau = -3$.  It is clear that these cannot all be congruent to each other $\pmod{4}$, hence perturbation of this type cannot yield perfect state transfer.

Conference graphs on $n$ vertices have $\theta,\tau = \frac{-1\pm\sqrt{n}}{2}$.  Either these are not integers, or when they are, it is clear they are not congruent $\pmod{4}$, hence perturbation of this type cannot yield perfect state transfer.
\qed

As an example of a strongly regular graph that does not belong to either of the infinite families given above, the complement of the Clebsch graph is a strongly regular graph with parameters $(16,10,6,6)$ and whose eigenvalues are \[k=10,\theta=2,\tau=-2,\] so Theorem \ref{thm:pstsrg-adj} implies that the perturbation on an edge can give perfect state transfer.

\section{Pretty Good State Transfer}\label{sec:pgst}

In this section, we show that in any primitive SRG $X$, there is a perturbation for which $X^{\beta,\gamma}$ has pgst between any pair of vertices.  Throughout this section, $X$ will be a strongly regular graph with parameters $(n,k,a,c)$ and eigenvalues $k,\theta,\tau$.  We will also assume that $X$ is primitive, that is, that $X$ is connected and not complete.  

The following lemma from \cite{BanchiCoutinhoGodsil} allows to study pretty good state transfer via a condition on the eigenvalues of the adjacency matrix.

\begin{lemma}\cite{BanchiCoutinhoGodsil}\label{lem:eig}
Let $u,v$ be vertices of $X$, and $A$ the adjacency matrix.  Then pretty good state transfer from $u$ to $v$ occurs at some time if and only if the following two conditions are satisfied:
\begin{enumerate}[(i)]
\item the vertices $u$ and $v$ are strongly cospectral; and 
\item letting $\{\lambda_i\}$ be the eigenvalues of $A$ corresponding to eigenvectors in the $U(-)$ space and $\{\mu_j\}$ are the eigenvalues for eigenvectors in $U(+)$,  if there exist integers $\ell_i$, $m_j$ such that if
\begin{align*}
\sum_i \ell_i\lambda_i +\sum_j m_j\mu_j = 0\\
\sum_i\ell_i +\sum_j m_j =0
\end{align*}
then
\[
\sum_i m_i \text{ is even.}
\]
\end{enumerate}
\end{lemma}

Following \cite{KemptonLippnerYauInvol}, for $X^{\beta,\gamma}$ where the perturbation is on adjacent vertices, let us define $P_+^{a}$ the polynomial whose roots are the eigenvalues corresponding to the $U(+)$ eigenspace for the, and $P_-^a$ the polynomial whose roots are the eigenvalues for the $U(-)$ eigenspace.  Let $P_+^n$ and $P_-^n$ denote the same polynomial in the non-adjacent version.

Let us denote
\begin{align*}
p_1(t) &= (t-k)(t-\theta)(t-\tau)\\
p_2(t) &= (t-\theta)(t-\tau)\\
q_{1,a}(t) &= t^2-(k+\theta+\tau-1)t+ k\theta+k\tau+k+2\theta\tau\\
q_{2,a}(t) &= t-(\theta+\tau+1)\\
q_{1,n}(t) &= t^2-(k+\theta+\tau)t+k\theta+k\tau+2k+2\theta\tau\\
q_{2,n}(t) &= t-(\theta+\tau)
\end{align*}
Then by Corollary \ref{lem:cubic-quadratic-equations},
 for adjacent vertices in a perturbation of a strongly regular graph
 \begin{align*}
 P_+^a(t) &= p_1(t) - (\beta+\gamma)q_{1,a}(t)\\
 P_-^a(t) &= p_2(t) - (\beta-\gamma)q_{2,a}(t)
 \end{align*} 
 and for non-adjacent vertices
 \begin{align*}
 P_+^n(t) &= p_1(t) - (\beta+\gamma)q_{1,n}(t)\\
 P_-^n(t) &= p_2(t) - (\beta-\gamma)q_{2,n}(t).
 \end{align*}

\begin{lemma}\label{lem:irred}
If $\beta+\gamma$ and $\beta-\gamma$ are both transcendental, then $P_+^a,P_-^a,P_+^n,P_-^n$ are all irreducible polynomials over $\Q(\beta,\gamma)$.
\end{lemma}
\proof
Since $\beta+\gamma$ and $\beta-\gamma$ are transcendental, and they appear as linear terms in the polynomials, if there is a factorization, there must be one rational factor. Thus we will be done if we can show that the pairs $p_1,q_{1,a}$, $p_2,q_{2,a}$, $p_1,q_{1,n}$, and $p_2,q_{2,n}$ are all relatively prime polynomials.  

If $p_1$ and $q_{1,a}$ are not relatively prime, then they share a root.  The roots of $p_1$ are $k,\theta,\tau$.  We have 
\begin{align*}
q_{1,a}(k) &= 2(k+\theta\tau) = 2c\\
q_{1,a}(\theta) &= (k+\theta)(\tau+1)\\
q_{1,a}(\tau) &= (k+\tau)(\theta+1).
\end{align*}
We know $c$ is not 0 for a connected strongly regular graph; for a connected strongly regular graph $k$ is strictly larger than $|\theta|$ or $|\tau|$, and  if $\theta$ or $\tau$ is -1, then that implies that the graph is complete.  Thus none of these is 0, so $p_1$ and $q_{1,a}$ are relatively prime.  

Similarly, the roots of $p_2$ are $\theta$ and $\tau$ and 
\begin{align*}
q_{2,a}(\theta) &= -(\tau+1)\\
q_{2,a}(\tau) &= -(\theta+1)
\end{align*}
and these are non-zero as above. So $p_2$ and $q_{2,a}$ are relatively prime.  

Similarly
\begin{align*}
q_{2,n}(\theta) = -\tau\\
q_{2,n}(\tau) = -\theta
\end{align*}
and neither $\theta$ not $\tau$ can be 0 in a connected strongly regular graph (this implies $k=c$ but $k>c$).  Therefore $p_2$ and $q_{2,n}$ are relatively prime.

Finally,
\begin{align*}
q_{2,n}(k) &= 2(k+\theta\tau) = 2c\\
q_{2,n}(\theta) &= k\tau+2k+\theta\tau=k\tau+c+k\\
q_{2,n}(\tau) &= k\theta+2k+\theta\tau=k\theta+c+k.
\end{align*}
As we saw above $c\neq0$.
\begin{claim}
For any strongly regular graph, $k\tau+c+k$ and $k\theta+c+k$ are non-zero.
\end{claim}
\proof
Suppose $k\tau+c+k = 0$.  Then $c=-k(\tau+1)$.  We know $c$ is a positive integer that is strictly less than $k$, so if $\tau\not\in\Q$ this is an immediate contradiction, but if $\tau\in\Q$, then $\tau$ is an integer, which also gives a contradiction. 

Similarly $k\theta+c+k$ is non-zero, giving the claim.
\qed

By the claim, the polynomials $p_2$ and $q_{2,n}$ are relatively prime, completing the proof.
\qed

\begin{theorem}\label{thm:pgst}
Given any pair of vertices $u$ and $v$ in any primitive strongly regular graph, then in $X^{\beta,\gamma}$, there is pretty good state transfer from $u$ to $v$ for any choice of $\beta,\gamma$ such that $P_+$ and $P_-$ are irreducible polynomials and $(\theta+\tau-(\beta-\gamma))/2\neq (k+\theta+\tau-(\beta+\gamma))/3$. In particular, the choice of $\beta$ and $\gamma$ that works is dense in the real numbers.
\end{theorem}
\proof
Choose $\beta$ and $\gamma$ satisfying the hypothesis of the theorem.  Let $\lambda_1,\lambda_2$ be the two roots of $P_-$ and $\mu_1,\mu_2,\mu_3$ the three roots of $P_+$.  
To use Lemma \ref{lem:eig}, we will examine the system
\begin{align*}
\ell_1\lambda_1+\ell_2\lambda_2+m_1\mu_1+m_2\mu_2+m_3\mu_3 &=0\\
\ell_1+\ell_2+m_1+m_2+m_3 &=0
\end{align*}
for $\ell_i,m_i$ integers.
We will treat the adjacent and non-adjacent cases together since the proof is the same.  We will suppress the $a$ or $n$ superscript and simply refer to the polynomials as $P_+$ and $P_-$.
In either the adjacent or non-adjacent case, the polynomials $P_+$ and $P_-$ are irreducible, $P_+$ with degree 3, and $P_-$ with degree 2.

We will apply techniques from \cite{KemptonLippnerYauInvol}.  
We will use a tool from Galois theory called the \emph{field trace} of a field extension.  For a field extension $K$ of $F$, we define $Tr_{K/F}: K \rightarrow F$ by 
\[
Tr_{K/F}(\alpha) = \sum_{g\in Gal(K/F)}g(\alpha).
\]  The field trace is the trace of the linear map taking $x\mapsto \alpha x$.  We record below a few basic facts about the field trace that will be useful.
\begin{enumerate}[(i)]
\item $Tr_{K/F}$ is a linear map.
\item For $\alpha\in F$, $Tr_{K/F}(\alpha) = [K:F]\alpha$.
\item For $K$ and extension of $L$, and extension of $F$, we have $Tr_{K/F} = Tr_{L/F}\circ Tr_{K/L}$.
\end{enumerate}

Let $F = \Q(\beta,\gamma)$.  Let $M/F$ be the splitting field for $P_+$, let $L/F$ be the splitting field for $P_-$, and let $K/F$ be the smallest field extension containing both $L$ and $M$.  By assumption, $P_+$ and $P_-$ are irreducible, so $L$ and $M$ are Galois extensions of $F$. Let us examine the field trace of the individual roots of $P_+$ and $P_-$.

Note that since our polynomials are irreducible, the Galois group $Gal(M/F)$ acts transitively on the roots, so taking $Tr_{M/F}(\mu_i)$, each of the other $\mu_j$'s will show up $|Gal(M/F)|/\mathrm{deg}~P_+$ times.  

Thus, for both the adjacent and non-adjacent cases, we have 
\[
Tr_{M/F}(\mu_i) = \frac{[M:F]}{3}(\mu_1+\mu_2+\mu_3) = \frac{[M:F]}{3}(k+\theta+\tau - (\beta+\gamma))
\]
and similarly
\[
Tr_{L/F}(\lambda_i) = \lambda_1+\lambda_2 = \theta+\tau-(\beta-\gamma).
\]

Now we have
\begin{align*}
0&=Tr_{K/F}\left(\ell_1\lambda_1+\ell_2\lambda_2+m_1\mu_1+m_2\mu_2+m_3\mu_3\right)\\
&=[K:L]Tr_{L/F}\left(\sum \ell_i\lambda_i\right)+[K:M]Tr_{M/F}\left(\sum m_j \mu_j\right)\\
&=[K:F]\left(\frac{\theta+\tau-(\beta-\gamma)}{2}\sum\ell_i + \frac{k+\theta+\tau-(\beta+\gamma)}{3}\sum m_j\right).
\end{align*}

Thus we have
\begin{align*}
\frac{\theta+\tau-(\beta-\gamma)}{2}\sum\ell_i + \frac{k+\theta+\tau-(\beta+\gamma)}{3}\sum m_j &=0\\
\sum\ell_i+\sum\mu_j &= 0.
\end{align*}
We have $\frac{\theta+\tau-(\beta-\gamma)}{2}\neq \frac{k+\theta+\tau-(\beta+\gamma)}{3}$, and thus these are two linearly independent equations in $\sum\ell_i$ and $\sum m_j$.  Thus 
\[
\sum\ell_i = \sum m_j = 0.
\]
In particular, each sum is even, so we get pretty good state transfer by Lemma \ref{lem:eig}.
The density claim follows from Lemma \ref{lem:irred}.
\qed

In Lemma \ref{lem:irred}, we assumed that $\beta$ and $\gamma$ were chosen so that $\beta+\gamma$ and $\beta-\gamma$ are both transcendental. The reason for this is that this is a generic choice for which we can prove that the polynomials $P_+$ and $P_-$ are irreducible, so that Theorem \ref{thm:pgst} applies.  However, there are certainly other, non-transcendental values of $\beta$ and $\gamma$ for which these polynomials are irreducible, but we do not know a nice way to characterize all of them.

For example, consider the Clebsch graph, which has parameters (16,5,0,2).  Its eigenvalues are $k=5, \theta=1,\tau=-3$. If we simply take consider a pair of non-adjacent vertices.  If we simply take $\beta=1,\gamma=0$, then the polynomials (\ref{eq:cubic-nonadj}) and (\ref{eq:quadratic-nonadj}) become
\begin{align*}
t^3-4t^2-10t+30&=0\\
t^2+t-5&=0
\end{align*}
which are irreducible polynomials over the rationals.  Therefore, our proof gives that there is pretty good state transfer between these vertices.  Since we took $\beta=1$, this corresponds to simply adding an edge between these vertices. So we have produced an example of a simple unweighted graph where pgst occurs.

\section{Conclusion and Open Problems}

In this paper, we characterized when the $(\beta,\gamma)$ perturbation on an edge $uv$, whose ends are cospectral vertices, of a graph results in strongly cospectral vertices. We specialize to strongly regular graphs, because they are a diverse class of graph where the eigenvalues and eigenspaces are well-understood. This would naturally extend to other families of graphs. In particular, Corollary \ref{cor:1wlkreg} gives that the eigenvalues of the $(\beta,\gamma)$ perturbation on $\{u,v\}$ of a 1-walk-regular graph $X$ is the same for any adjacent pair of vertices $u,v$. This is not apparent a priori, but follows as a consequence of the analysis. One can ask for other classes of graph does this property hold. 

Furthermore, Theorems \ref{thm:pstsrg-nonadj} and \ref{thm:pstsrg-adj} characterize when a perturbation of a strongly regular graph can yield perfect state transfer in the particular case $\beta=-\gamma$.  This condition on $\beta$ and $\gamma$ simplified the analysis, but it is still of interest to ask if a perturbation with some more general choice of $\beta$ and $\gamma$ could yield perfect state transfer in other strongly regular graphs. 

Of particular interest along these lines is the question of whether there are strongly regular graphs where, taking $\beta=\pm1$ and $\gamma=0$ yields perfect state transfer.  This corresponds to addition or deletion of an edge, and would thus yield unweighted graphs with no vertex potential on which there is perfect state transfer.

A natural extension would be to ask for similar result about the perturbations of hypercubes. Hypercubes of order $n$ are not strongly regular for $n>2$, but are $1$-walk-regular.  

\end{document}